\documentclass[11pt]{article}
\usepackage{amsfonts}

\newcommand{\oper}[2]{\newcommand{#1}{\mathop{\kern0pt\mathrm{#2}}\nolimits} }
\oper{\tr}{tr}
\oper{\adj}{adj}
\oper{\Div}{div}
\oper{\ad}{ad}
\oper{\Ad}{Ad}
\oper{\End}{End}
\oper{\Hom}{Hom}
\oper{\Aut}{Aut}
\oper{\SO}{SO}
\oper{\SP}{Sp}
\oper{\SU}{SU}
\oper{\GL}{GL}
\oper{\id}{I}
\oper{\ext}{Ext}
\oper{\rank}{rank}
\oper{\diag}{Diag}
\oper{\sgn}{sgn}

\newcommand{\ka}{K\"ahler}
\newcommand{\norm}[1]{\Vert #1\Vert}  
\newcommand{\abs}[1]{\vert #1\vert} 
\newcommand{\h}[1]{\hat #1} 
\newcommand{\lie}[1]{\mathfrak{#1}}

\newcommand{\htheta}{{\hat\theta}}
\newcommand{\eqref}[1]{(\ref{#1})}

\newcommand{\ov}{\overline}

\newcommand{\del}{\partial}
\newcommand{\delb}{\overline{\partial}}

\newcommand{\CC}{\mathbb{C}}   
\newcommand{\HH}{\mathbb{H}}   
\newcommand{\RR}{\mathbb{R}}

\newtheorem{theorem}{Theorem}
\newtheorem{corollary}[theorem]{Corollary}
\newtheorem{proposition}[theorem]{Proposition}
\newtheorem{definition}[theorem]{Definition}
\newtheorem{lemma}[theorem]{Lemma}

\newcommand{\bproof}{\noindent{\it Proof: }}
\newcommand{\eproof}{\  q.~e.~d. \vspace{0.2in}}

\begin{document}

\title{Potential One-Forms for Hyperk\"ahler Structures\\
with Torsion}
\author{
Liviu Ornea 
\thanks{Member of \textsc{Edge}, Research Training Network
\textsc{hprn-ct-\oldstylenums{2000}-\oldstylenums{00101}}, supported by
The European Human Potential Programme.
Address: University of Bucharest, Faculty of Mathematics
14 Academiei str., 
70109 Bucharest, Romania.
E-mail: lornea@imar.ro}
\and 
Yat Sun Poon
\thanks{Partially supported by NSF DMS-0204002. 
Address: Department of Mathematics, University of California at Riverside, Riverside, CA 92521, U.S.A..
E-mail: ypoon@math.ucr.edu}
\and
Andrew Swann
\thanks{Member of \textsc{Edge}, Research Training Network
\textsc{hprn-ct-\oldstylenums{2000}-\oldstylenums{00101}}, supported by
The European Human Potential Programme.
Address: Department of Mathematics and Computer Science,
University of Southern Denmark,
Campusvej 55,
DK-5230 Odense M,
Denmark.
E-mail: swann@imada.sdu.dk}
}
\date{November 22, 2002}
\maketitle

\noindent{\bf Abstract}\hspace{.15in}
It is shown that an HKT-space with closed parallel potential 1-form has
$D(2,1;-1)$-symmetry.  Every locally conformally hyperk\"ahler manifold
generates this type of geometry.  The HKT-spaces with closed parallel
potential 1-form arising in this way are characterized by their symmetries
and an inhomogeneous cubic condition on their torsion.

\section*{Introduction}
 HKT-geometry is a metric geometry with multiple complex structures that
 arises in various physical theories, including supersymmetric non-linear
 sigma models, type IIA string theory, and black hole moduli.  Good references for the
 physical background are \cite{GPS} and \cite{MS} and the citations
 therein.  For a mathematical approach, we refer the reader to \cite{GP}.
Since the geometry is typically hyperhermitian and non-K\"ahlerian,
 it is of great interest and challenging
to find potential functions \cite{MS}.

 In the context of multi-particle quantum mechanics, Michelson and
 Strominger studied the phenomenon of superconformal symmetry.  Motivated
 by application to dynamics of black holes \cite{MS2},
 they demonstrated in \cite{MS} a relation between a $D(2,1;\alpha)$
 superconformal symmetry and classical differential geometry on
 HKT-manifolds.
Given supersymmetry such as this, potential functions are already found
\cite{PS:a} \cite{PS:b}.

On the other hand, a maximum principle argument shows that potential
functions could not exist on compact manifolds \cite{GP}. We therefore
replace locally defined potential functions by a globally defined closed
1-form in our consideration (see Definition \ref{potential form}). We focus
on the case when the potential 1-form is parallel with respect to the
HKT-connection in this investigation. 
Combining Corollary \ref{main corollary} and Proposition \ref{converse}, we obtain 
the following result in this direction. 

{\it If $V$ is the dual vector field of a closed parallel potential 1-form
$\theta$ of the HKT-space with metric $\hat{g}$ and hypercomplex structures
$I_1, I_2, I_3$, then
\[
d\theta=0,
\quad
{\mathcal L}_V{\hat{g}}=0,
\quad
{\mathcal L}_{I_rV}{\hat{g}}=0,
\quad 
{\mathcal L}_{I_rV}I_s=\epsilon^{rst}I_t.
\]
Conversely, if there is such a vector field on an HKT-space, then the dual 1-form is a parallel potential function.}

Due to a theorem of Michelson and Strominger \cite{MS}, this type of
symmetry is a degenerate version of $D(2,1;\alpha)$ symmetry, namely
$D(2,1;-1)$. Since the above symmetry makes sense on the HKT-space, we
shall refer to it as $D(2,1;-1)$-symmetry in this paper
 despite an apparent singularity that occurs in the 
structural equations \cite[3.44]{MS}. Due to an isomorphism among superalgebras
with different parameters
\cite[Proposition 2.5.4]{Kac}, $D(2,1;-1)$ is isomorphic to $D(2,1;0)$ and $D(2,1;\infty)$. 
This equivalent class of superalgebras is featured to have one decoupled $SU(2)$.
In this paper, we interpret $D(2,1;-1)$ symmetry  after  Michelson and
Strominger's theorem \cite[3.56]{MS}.  
A precise description is given in Definition
\ref{def-d21}.  Through a construction, we shall prove the following
observation.

\it If $(M,g, I_1, I_2, I_3)$ is a locally conformally hyperk\"ahler
manifold whose Lee form is parallel with respect to the Levi-Civita
connection, then there exists an HKT-metric $\hat{g}$ such that the Lee form
of $g$ is a potential 1-form of $\hat{g}$, and is parallel with respect to
the HKT-connection of $\hat{g}$.  \rm

As a result in potential theory, the above observation supplements what is
already known for HKT-spaces with $D(2, 1;\alpha)$-symmetry when $\alpha\neq
-1, 0, \infty$. From a geometric perspective, it implicitly links
HKT-geometry to Weyl geometry, quarternionic geometry and Sasakian geometry
through the theory of locally conformally hyperk\"ahler manifolds.

We conclude with a discussion on how to distinguish the class of HKT-spaces
associated to locally conformally hyperk\"ahler manifolds.

 Throughout this article
we adopt the conventions in \cite{Besse} and
\cite{Gauduchon2}.  Here we warn casual readers that the concerned metrics
for locally conformally hyperk\"ahler structure and its associated
HKT-structure are in different conformal classes.

\section{HKT-Manifolds}
A Hermitian structure on a smooth manifold $M$ consists of a Riemannian metric $\hat{g}$ and an integrable complex structure $J$ such that for any tangent vectors $X$ and $Y$ on the manifold $M$,
\[
{\hat{g}}(JX, JY)={\hat{g}}(X,Y).
\] 
A triple of integrable
complex structure $I_r$, $r=1,2,3$, forms a hypercomplex structure on the manifold $M$ if
they satisfy the quaternion relations:
\[
I_1^2=I_2^2=I_3^2=\id, \quad I_1I_2=I_3=-I_2I_1.
\]

If each complex structure $I_r$ with the metric $\hat{g}$ forms a Hermitian structure, then
$(M,{\hat{g}},I_1,I_2,I_3)$ is said to be a hyperhermitian manifold. 

We denote ${\hat{F}}_r$ the fundamental two-form associated to the complex structure $I_r$ and we observe the convention:
\[
{\hat{F}}_r(X,Y)={\hat{g}}(I_rX,Y).
\]
For a $k$-form $\omega$ let 
\begin{equation}\label{zero}
(I_r\omega)(X_1,...,X_k)=(-1)^k\omega(I_rX_1,...,I_rX_k).
\end{equation}
The complex operators $d_r$, $\del_r$ and $\delb_r$ are respectively defined as:
\[
 d_r\omega =(-1)^kI_rdI_r\omega\quad \mbox{for a $k$-form $\omega$,}
\quad \del_r = \frac{1}{2}(d+id_r), \quad \delb_r=\frac{1}{2}(d-id_r).
\]

\begin{definition}
A linear connection $D$ with torsion tensor $T^D$ on $M$ is called 
\emph{hyperk\"ahler with torsion} if 
\par (i) it is hyperhermitian: $DI_1=DI_2=DI_3=0$, $Dg=0$ and
\par (ii) the tensor field $c$ defined by
$c(X,Y,Z)={\hat g}(T^D(X,Y),Z)$ is a $3$-form.
\end{definition}

Such a connection is denoted HKT by physicists  \cite{GPS} \cite{MS}
and we shall preserve this name. 
 Among mathematicians, HKT-connection is also known as Bismut connection for each of the
complex structures $I_r$  \cite{Gauduchon2}. 
Using the characterization of the Bismut connection and the
fact that  it is  uniquely associated to a Hermitian
structure, one obtains the following equivalent observation \cite{GP} \cite{GPS}:
\begin{proposition}\label{df}
On any hyperhermitian manifold $(M,{\hat{g}},I_1,I_2,I_3)$, the following two conditions are
equivalent 
\par (i) $d_1{\hat F}_1=d_2{\hat F}_2=d_3{\hat F}_3.$
\par (ii) $\del_1({\hat F}_2+i{\hat F}_3)=0.$

\noindent An HKT-connection exists if and only if one of the above two
conditions is satisfied.  When it exists, it is unique.
\end{proposition}

As demonstrated in \cite{MS}, an efficient way for constructing examples of HKT structures is the use of
HKT potentials. 
These are generalizations of hyperk\"ahler potentials \cite{GP}.
\begin{definition}
  Let $(M, {\hat{g}},I_1,I_2,I_3)$ be an HKT manifold. A (possibly locally
  defined) function $\mu:U\subseteq M\rightarrow \RR$ is a \emph{potential
  function} for the HKT structure if
\begin{equation}\label{pot}
{\hat F}_1=\frac{1}{2}(dd_1+d_2d_3)\mu, \quad {\hat F}_2=\frac{1}{2}(dd_2+d_3d_1)\mu, \quad 
{\hat F}_3=\frac{1}{2}(dd_3+d_1d_2)\mu.
\end{equation}
\end{definition}
Alternatively, the potential function $\mu$ is characterized by 
\begin{equation}\label{pott}
{\hat F}_2+i{\hat F}_3=2\del_1I_2\delb_1\mu.
\end{equation}

Potential functions do not always exist. When one exists, the torsion form
of an HKT structure deriving from a potential $\mu$ is:
$$c=-\frac{1}{2}d_1d_2d_3\mu=-d_1{\hat F}_1=-d_2{\hat F}_2=-d_3{\hat
F}_3.$$
As an example, the function $\log \sum_i\abs{z_i}^2$ is an HKT
potential on $\CC^{2n}\backslash\{0\}$. Moreover, it descends locally to
the Hopf manifold $S^1\times S^{4n-1}$.

This should be noted that like K\"ahler potentials, HKT-potentials could
not exist globally on compact manifolds due to a typical maximum principle
argument \cite{GP}. Moreover, a generic HKT-manifold is non-K\"ahlerian and
the $\partial{\overline \partial}$-lemma is not applicable. Therefore, we
propose to develop a global version of potential theory through the
Poincar\'e Lemma for 1-forms.
\begin{definition}\label{potential form} A one-form $\omega$ is a potential
  1-form for an HKT-manifold $(M, {\hat{g}},I_1,I_2,I_3)$ if the
  fundamental two-forms are given by
\begin{equation}
  {\hat F}_1=\frac{1}{2}(d\omega_1+d_2\omega_3), \quad {\hat
  F}_2=\frac{1}{2}(d\omega_2+d_3\omega_1), \quad  
  {\hat F}_3=\frac{1}{2}(d\omega_3+d_1\omega_2),
\end{equation}
where $\omega_r:=I_r\omega$. A potential 1-form is closed if $d\omega=0$.
\end{definition}

In such terminology, the HKT-structure on Hopf manifolds has a globally
defined potential 1-form.  Implicitly, 
 Poincar\'e
Lemma provides the locally
defined potential functions whenever a potential 1-form exists and is
closed.  Moreover, the torsion 3-form is now given by
\begin{equation}
c=-\frac12d_1d_2\omega_3=-\frac12d_2d_3\omega_1=-\frac12d_3d_1\omega_2.
\end{equation}

\section{Parallel Potential Forms}
In this section, we analyze the structure of HKT-spaces with parallel
potential 1-forms.  Since HKT-connections are Riemannian connections,
vector fields dual to parallel potential forms are parallel. Therefore, we
extend our investigation to parallel vector fields in general briefly,
before we focus again on potential 1-forms and their dual vector fields.

\begin{lemma}\label{parallelism}
  Let $V$ be a vector field on an HKT-space. The following statements are
  equivalent:
\par (i) $V$ is parallel with respect to the HKT-connection $D$.
\par (ii) $V, I_1V, I_2V, I_3V$ are parallel with respect to the HKT-connection $D$.
\par (iii) $V, I_1V, I_2V, I_3V$ are Killing vector fields with respect to the
HKT-metric.
\end{lemma}
\bproof
Since HKT-connection preserves the hypercomplex structure, the equivalence between
the first two statements are obvious.

For any vector fields $W, Y, Z$, as $D$ is a metric connection, we have the identity
\begin{eqnarray*}
{\mathcal L}_W{\hat g}(Y, Z)
&=& {\hat g}(D_YW, Z)+{\hat g}(Y, D_ZW)+
{\hat g}(T^D(W, Y), Z)+{\hat g}(Y, T^D(W, Z))\\
&=& {\hat g}(D_YW, Z)+{\hat g}(Y, D_ZW)+
c(W, Y, Z)+c(Y, W, Z).
\end{eqnarray*}
Since $c$ is totally skew, we have 
\begin{equation}
{\mathcal L}_W{\hat g}(Y, Z)
= {\hat g}(D_YW, Z)+{\hat g}(Y, D_ZW)
\end{equation}
Applying this identity to the vector fields $V, I_1V, I_2V, I_3V$ and 
using the fact that the HKT-connection preserves the hypercomplex structure, 
we derive the implication from the second statement to the third.

Conversely, if the vector fields $V, I_1V, I_2V, I_3V$ are Killing, we apply the above
identity to $V$ to conclude that the symmetric part of $D V$ is equal to zero.
Let $\beta$ be the skew-symmetric part of 
 $D V$, i.e., $D V=\beta$.
Since the connection preserves the complex structures, the above identity
is equivalent to
\begin{equation}
{\hat g}(D_Y(I_rV), Z)={\hat g}(I_rD_YV, Z)=-\beta(Y, I_rZ).
\end{equation}
On the other hand, as the vector fields are Killing, 
\begin{equation}
{\hat g}(D_Y(I_rV), Z)+{\hat g}(D_Z(I_rV), Y)=({\mathcal L}_{(I_rV)}{\hat g})(Y, Z)=0.
\end{equation}
Therefore, $\beta(Y, I_rZ)+\beta(Z, I_rY)=0.$
Then
\begin{eqnarray*}
\beta(Y, I_1Z) &=& -\beta(Z, I_1Y)=\beta(I_1Y, Z)=\beta(I_2I_3Y, Z)\\
&=& \beta(I_3Y, I_2Z)=\beta(Y, I_3I_2Z)=-\beta(Y, I_1Z).
\end{eqnarray*}
Therefore, $\beta=0$. This implies that 
$DV=0$.
\eproof

\begin{lemma}\label{diff} Suppose that $V$ is a parallel vector field with respect to the HKT-connection
$D$. Let $\hat\theta$ be its dual 1-form with respect to $\hat g$. Then
\begin{equation}
d{\hat\theta}=\iota_Vc,
\qquad
d{\hat\theta}_r=\iota_{I_rV}c.
\end{equation}
\end{lemma}
\bproof
Let $0\leq m\leq 3$. Let $I_0$ denote the identity endomorphism on tangent space.
For any vector fields $X$ and $Y$, 
\begin{eqnarray*}
d{\hat\theta}_m(X, Y)&=& X({\hat\theta}_m(Y))-Y({\hat\theta}_m(X))-{\hat\theta}_m([X, Y])\\
&=& X({\hat g}(I_mV, Y))-Y({\hat g}(I_mV, X))-g(I_mV, [X, Y])\\
&=& {\hat g}(I_mV, D_XY)-{\hat g}(I_mV, D_YX)-{\hat g}(I_mV, [X, Y])\\
&=& {\hat g}(I_mV, T^D (X, Y))=c(I_mV, X, Y).
\end{eqnarray*}
\eproof

\begin{lemma}\label{contraction} Suppose that $V$ is a parallel vector field with respect to the HKT-connection
$D$. It is parallel with respect to the Levi-Civita connection $\hat\nabla$ of the metric 
$\hat g$ if and only if $\iota_Vc=0$.
\end{lemma}
\bproof
This is due to the identity
${\hat g}({\hat\nabla}_XV, Y)={\hat g}(D_XV, Y)+c(X, V, Y)=c(X, V, Y)$.

\

Next we investigate the behavior of the vector fields $V, I_1V, I_2V, I_3V$ with respect to the hypercomplex
structure $\{I_1, I_2, I_3\}$.

\begin{lemma}\label{rotating} If $-2\hat\theta$ is a closed 
potential 1-form and is parallel with respect to the HKT-connection, then 
${\mathcal L}_VI_r=0$ and  ${\mathcal L}_{I_rV}I_s=\epsilon^{rst}I_t$.
\end{lemma}
\bproof Since the vector fields $V, I_1V, I_2V, I_3V$ are Killing vector fields, it suffices to show that
${\mathcal L}_V{\hat F}_r=0,$ and ${\mathcal L}_{I_rV}{\hat F}_s=\epsilon^{rst}{\hat F}_t.$

In the following computation, we use the results in Lemma \ref{diff} extensively. 
For any tangent vectors $X$ and $Y$, 
\[
\iota_V d{\hat F}_r(X, Y)= (\iota_VI_rc)(X,Y)
= -c(I_rV, I_rX, I_rY)=-d{\hat\theta}_r(I_rX, I_rY)=-I_rd{\hat\theta}_r(X,Y).
\]
On the other hand, 
$
\iota_V{\hat F}_r(X)={\hat g}(I_rV, X)={\hat\theta}_r(X).
$
Therefore,
\[
{\mathcal L}_V{\hat F}_r
= \iota_V d{\hat F}_r+d\iota_V{\hat F}_r=-I_rd{\hat\theta}_r+d{\hat\theta}_r.
\]
As the torsion form is of type $(1,2)+(2,1)$ with respect to all $I_r$, 
\begin{equation}\label{type}
c(Z, X, Y)=c(Z, I_rX, I_rY)+c(I_rZ, X, I_rY)+c(I_rZ, I_rX, Y).
\end{equation}
Substitute $Z$ by $I_rV$  and apply Lemma \ref{diff}, we have 
\[
d{\htheta}_r(X, Y)=I_rd{\htheta}_r(X,Y)-d\htheta(X, I_rY)-d\htheta(I_rX, Y).
\]
Therefore, 
$
{\mathcal L}_V{\hat F}_r(X,Y)=-d\htheta(X, I_rY)-d\htheta(I_rX, Y).
$
As $\htheta$ is closed,  
${\mathcal L}_V{\hat F}_r=0$.
Next, 
\begin{equation}\label{I-F}
\iota_{I_rV}{\hat F}_r(X)={\hat F}_r(I_rV, X)={\hat g}(I^2_rV, X)=-\htheta(X).
\end{equation}
With Lemma \ref{diff}, we have
\begin{equation}
\iota_{I_rV}d{\hat F}_r(X,Y) = \iota_{I_rV}I_rc(X, Y)
= -c(I^2_rV, I_rX, I_rY)=c(V, I_rX, I_rY)=I_rd\htheta (X,Y).
\label{Id-theta}
\end{equation}
Therefore, 
\begin{equation}\label{Lie-F}
{\mathcal L}_{I_rV}{\hat F}_r=\iota_{I_rV}d{\hat F}_r+d\iota_{I_rV}{\hat F}_r
=I_rd\htheta-d\htheta.
\end{equation}
Since $d\htheta=0$, ${\mathcal L}_{I_rV}{\hat F}_r=0$.
Finally,  
\begin{equation}\label{I1-F2}
\iota_{I_1V}{\hat F}_2(X)={\hat F}_2(I_1V, X)={\hat g}(I_2I_1V, X)=-\htheta_3(X).
\end{equation}
By Lemma \ref{diff} and (\ref{type}),
\begin{eqnarray}
 \iota_{I_1V}d{\hat F}_2(X,Y)&=& \iota_{I_1V}I_2c(X, Y)=I_2c(I_1V, X,Y)
 =c(I_3V, I_2X, I_2Y)
\nonumber\\
&=& c(I_3V, I_3I_2X, I_3I_2Y)+c(I^2_3V, I_2X, I_3I_2Y)+c(I^2_3V, I_3I_2X, I_2Y)
\nonumber\\
&=& c(I_3V, I_1X, I_1Y)+c(V, I_2X, I_1Y)+c(V, I_1X, I_2Y)
\nonumber\\
&=& I_1d\htheta_3(X, Y)+d\htheta (I_2X, I_1Y)+d\htheta (I_1X, I_2Y). 
\label{iota-1-F2}
\end{eqnarray}
Therefore, 
\begin{equation}\label{Lie-1-F2}
{\mathcal L}_{I_1V}{\hat F}_2(X,Y)=-d\htheta_3(X,Y)+I_1d\htheta_3(X,Y)+d\htheta (I_2X, I_1Y)+d\htheta (I_1X, I_2Y).
\end{equation}
On the other hand, if $-2\htheta$ is a potential 1-form, then $d\htheta=0$. It follows that
\[
{\mathcal L}_{I_1V}{\hat F}_2=-d\htheta_3+I_1d\htheta_3.
\]
In addition,
\[
{\hat F}_3=\frac12 (d(-2\htheta_3)+d_1(-2\htheta_2))
=-d\htheta_3+I_1dI_1I_2\htheta
=-d\htheta_3+I_1d\htheta_3.
\]
Therefore, ${\mathcal L}_{I_1V}{\hat F}_2={\hat F}_3$.
\eproof

Summarizing the results in Lemma \ref{parallelism} and Lemma \ref{rotating} in the context of parallel
potential 1-forms, we have the next result. 
\begin{corollary}\label{main corollary} Suppose that $-2\hat\theta$ is a closed potential 1-form and parallel with respect to
the HKT-connection. If $V$ is the dual of $\htheta$ with respect to the HKT-metric $\hat{g}$, then 
\begin{equation}\label{d21}
{\mathcal L}_V{\hat{g}}=0,
\quad
{\mathcal L}_{I_rV}{\hat{g}}=0,
\quad 
{\mathcal L}_{I_rV}I_s=\epsilon^{rst}I_t.
\end{equation}
\end{corollary}

Comparing with \cite[(3.56)]{MS} and keeping in mind that the dual 1-form $\htheta$ is closed,
we conclude that the HKT-space in
question is induced by the $D(2,1;-1)$-supersymmetry. Although such supersymmetry is singular as
seen in \cite[(3.44)]{MS}, we retain the notion of $D(2,1;-1)$-symmetry. To be precise, we make a definition.

\begin{definition}\label{def-d21} A $D(2,1;-1)$-symmetry on an HKT-space is
  a vector field $V$ satisfying the conditions in (\ref{d21}) and whose
  dual 1-form $\htheta$ is closed.
\end{definition}

In previous investigation on potential functions \cite{PS:a} \cite{PS:b}, such symmetry was not 
 extensively studied
due to degeneracy of supersymmetry. Below is a remedy. 

\begin{proposition}\label{converse} Suppose that a vector field $V$
  generates a $D(2,1;-1)$-symmetry on an HKT-space. Let $\hat\theta$ be the
  dual vector field. Then $-2\hat\theta$ is a parallel potential 1-form. In
  particular, local potential function exists.
\end{proposition}
\bproof
By definition, $V, I_1V, I_2V, I_3V$ are Killing vector fields. By Lemma \ref{parallelism},
$V$ is parallel with respect to the HKT-connection. In particular, Lemma \ref{diff} is 
applicable. With it, we obtain equation (\ref{iota-1-F2}). With identity (\ref{I1-F2}), we
obtain equation (\ref{Lie-1-F2}). Since $\htheta$ is closed, 
${\mathcal L}_{I_1V}{\hat F}_2=-d\htheta_3+I_1d\htheta_3$. 
On the other hand, as $I_1V$ is a Killing vector field and 
${\mathcal L}_{I_1V}I_2=I_3$, it follows that ${\mathcal L}_{I_1V}{\hat{F}}_2={\hat{F}}_3$. Therefore, 
\[
{\hat{F}}_3=-d\htheta_3+I_1d\htheta_3
=\frac12 (d (-2\htheta_3)+d_1(-2\htheta_2)).
\]
The above calculation is repeated with the indices permuted to conclude that $-2\htheta$ is a
potential 1-form.
\eproof

\noindent{\bf Remark:} 
  By Lemma \ref{diff} and Lemma \ref{contraction}, closedness of $\htheta$
  along with the parallelism of the dual vector field $V$ together implies
  the vector field of symmetry is parallel with respect to the Levi-Civita
  connection of the HKT-metric $\hat{g}$. In view of Lemma \ref{rotating},
  it implies that ${\mathcal L}_VI_r=0$.

\section{Locally Conformally Hyperk\"ahler Manifolds}

Locally conformally hyperk\"ahler manifolds have been studied in relation
to Weyl geometry, quaternionic geometry as well as Sasakian geometry
\cite{OP1} \cite{PPS1}.  In this section, we demonstrate a way to generate
HKT-structures with $D(2,1;-1)$-symmetry and parallel potential 1-form from
a locally conformally hyperk\"ahler structure.  We begin our investigation
with a review of definitions.

\begin{definition}
\par (i) A hyperhermitian manifold $(M,g,I_1,I_2,I_3)$ is cal\-led 
\emph{hyperk\"ahler} if the Levi-Civita connection of $g$ parallelizes each 
complex structure $I_r$: $\nabla I_r=0$.
\par (ii) A hyperhermitian manifold $(M,g,I_1,I_2,I_3)$ is called
\emph{locally conformally hyperk\"ahler} if there exists an open 
cover $\{U_i\}$ such that the restriction of the metric to each $U_i$ is 
conformal to a local hyperk\"ahler metric $g_i$: 
\begin{equation}\label{unu}
g|_{U_i}=e^{f_i}g_i, \quad f_i\in \mathcal{C}^\infty{U_i}.
\end{equation}
\end{definition}

\

 We shall focus on the second notion. Taking $\theta|_{U_i}=df_i$, the condition 
(\ref{unu}) is equivalent to the existence 
 of a globally defined one-form $\theta$ satisfying the integrability conditions:
\begin{equation}\label{dFr}
dF_r=\theta\wedge F_r, \quad r=1,2,3.
\end{equation}

        The standard example of locally conformally hyperk\"ahler manifold is the Hopf manifold 
$H^n_\HH=(\HH\backslash\{0\})/\Gamma_2$, where 
$\Gamma_2$ is the cyclic group generated by the quaternionic automorphism 
$(q_1,...,q_n)\mapsto (2q_1,...,2q_n)$. The hypercomplex structure 
of $\HH^n$ is easily seen to descend to $H_\HH^n$. Moreover, 
the globally conformal hyperk\"ahler metric $(\sum_iq_i\ov{q}_i)^{-1}
\sum_idq_i\otimes d\ov{q}_i$ on $\HH^n\backslash\{0\}$ is invariant to the action of 
$\Gamma_2$, hence induces a locally conformally hyperk\"ahler metric on 
the Hopf manifold with Lee form 
$$
\theta=-\frac{\sum_i(q_id\ov{q}_i+\ov{q}_idq_i)}{\sum_iq_i\ov{q}_i}.
$$
Note that,  as in the 
complex case, $H^n_\HH$ is diffeomorphic with a product of spheres 
$S^1\times S^{4n-1}$. Consequently, 
its first Betti number is $1$ and it cannot admit any hyperk\"ahler metric.
        Other examples are presented in \cite{OP1} where also a complete 
classification of compact homogeneous locally conformally hyperk\"ahler manifolds is
given.

One should note that locally conformally hyperk\"ahler manifolds are
hyperhermitian Weyl and as such, Einstein-Weyl Ricci-flat (here, the
conformal class is that of $g$ and the Weyl connection is constructed out
of the Levi-Civita connection of $g$ and the Lee form). Hence, if compact,
one applies a well-known result of Gauduchon \cite{Gauduchon1}, to obtain
the existence of a metric $g_0$, conformal with $g$ and having the Lee form
parallel with respect to the Levi-Civita connection of $g_0$. The metric we
just wrote on the Hopf manifold has this property.  Therefore, when working
with compact locally conformally hyperk\"ahler manifolds, one can always
assume the metric with parallel Lee form.  We shall need the following
computational result \cite{OP1}:

\begin{lemma}
Let $(M,g,I_1,I_2,I_3)$ be a locally conformally hyperk\"ahler manifold with 
parallel Lee form $\theta$. Let $\theta_r=I_r\theta$. Assume that $\theta$ has unit length. Then 
\begin{equation}\label{trei}
d\theta_r=\theta\wedge\theta_r-F_r.
\end{equation}
\end{lemma}

It should be noted that
 the unit length condition may achieved by rescaling $g$ by a homothety
 and that
\begin{equation}
I_rd\theta_r=I_r\theta\wedge I_r\theta_r-I_rF_r=-\theta_r\wedge\theta-F_r
=d\theta_r.
\end{equation}
Also,
\begin{equation}
I_rdF_r=I_r\theta\wedge I_rF_r=\theta_r\wedge F_r.
\end{equation}

That the Hopf manifolds admit HKT structures is not by chance. We can state:
\begin{theorem}
Let $(M,g,I_1,I_2,I_3)$ be a locally conformally hyper\-k\"a\-hler manifold with parallel
Lee form $\theta$. Assume that $\theta$ has unit length. Then the metric 
\begin{equation}\label{main}
\h{g}=g-\frac{1}{2}\{\theta\otimes\theta+
\theta_1\otimes\theta_1+\theta_2\otimes\theta_2+\theta_3\otimes\theta_3\}
\end{equation}
is HKT. Moreover, $\theta$ is a closed potential 1-form for $\h{g}$.
\end{theorem}
\bproof Let $g_2= \theta\otimes\theta+
\theta_1\otimes\theta_1+\theta_2\otimes\theta_2+\theta_3\otimes\theta_3$ 
be the restriction of the metric $g$ on the quaternionic span of the vector field $V$. 
Let $g_1$ be the restriction of the metric $g$ on the orthogonal complement of 
the quaternionic span of $V$. Then
the metric $g$ pointwisely and smoothly splits into two parts
$g=g_1+g_2.$ 
Since the norm of $\theta$ and its dual vector field $V$ have unit length with respect to $g$, 
the bilinear form $\h{g}$ is equal to $g_1+\frac{1}{2}g_2$. In particular,
this is a Riemannian metric.

Note first that, due to
(\ref{zero}) we have:
\begin{equation}
I_rF_r =F_r, \quad I_rF_s=-F_s \quad \mbox{for $r\neq s$},
\quad
I_r\theta_s =\epsilon^{rst}\theta_t.
\end{equation}
As a matter of convention, for exterior products we use that 
\begin{equation}
\alpha_1\wedge\cdots\wedge\alpha_n(X_1, \cdots, X_n):=\det (\alpha_i(X_j)).
\end{equation}
In particular,
$\theta\wedge\theta_1=\theta\otimes\theta_1-\theta_1\otimes\theta.$
From the definitions and (\ref{main}), 
\begin{equation}\label{eq:hat-F}
\h{F}_1=F_1-\frac{1}{2}\{\theta\wedge\theta_1+\theta_2\wedge\theta_3\}.
\end{equation}
Now we have successively, using $d\theta=0$, $dF_r=\theta\wedge F_r$ and formula (\ref{trei}):
\begin{eqnarray}
d\h{F}_1&=&dF_1-\frac{1}{2}\{d\theta\wedge\theta_1-\theta\wedge d\theta_1+d\theta_2\wedge\theta_3-\theta_2\wedge d\theta_3\}\nonumber\\
&=&dF_1-\frac{1}{2}\{-\theta\wedge(\theta\wedge\theta_1-F_1)+(\theta\wedge\theta_2-F_2)\wedge\theta_3 -\theta_2\wedge(\theta\wedge\theta_3-F_3)\}\nonumber\\
&=& \frac{1}{2}\{\theta\wedge F_1-2\theta\wedge\theta_2
   \wedge\theta_3+\theta_3\wedge F_2-\theta_2\wedge F_3\}. \label{dF-hat}\\
I_1d\h{F}_1&=&\frac{1}{2}\{\theta_1\wedge I_1F_1 -2\theta_1\wedge I_1\theta_2\wedge I_1\theta_3 +I_1\theta_3\wedge I_1F_2-I_1\theta_2\wedge I_1F_3\}\nonumber\\
&=&\frac{1}{2}\{\theta_1\wedge F_1+\theta_2\wedge F_2+\theta_3\wedge F_3-2\theta_1\wedge \theta_2\wedge\theta_3\}. \label{eq:c-theta-F}
\end{eqnarray}
The above formula is symmetric in the indices $1,2,3$. Due to Proposition
\ref{df}, $\hat{g}$ is an HKT-metric.

We prove the assertion on potential one-form by demonstrating that
any locally defined function $f$ with $df=\theta$ is a potential function. 
\begin{eqnarray*}
  \delb_1f &=& \frac{1}{2}(df-iI_1df) = \frac{1}{2}(\theta-iI_1\theta) =
  \frac{1}{2}(\theta-i\theta_1),\\   
  I_2\delb_1f &=& \frac{1}{2}(I_2\theta-iI_2\theta_1) =
  \frac{1}{2}(\theta_2+i\theta_3),\\ 
  \del_1I_2\delb_1f &=&
  \frac{1}{4}(d\theta_2+id\theta_3-iI_1d(I_1\theta_2+iI_1\theta_3))
=
  \frac{1}{4}(d\theta_2+id\theta_3-iI_1d(\theta_3-i\theta_2))
\\
  &=& \frac{1}{4}(\theta\wedge\theta_2-F_2+i(\theta\wedge\theta_3-F_3)
  - iI_1(\theta\wedge\theta_3-F_3)-I_1(\theta\wedge\theta_2-F_2))\\
  &=& -\frac{1}{2}(F_2+iF_3) + \frac{1}{4}(\theta+i\theta_1) \wedge
  (\theta_2+i\theta_3). 
\end{eqnarray*}
On the other hand, 
$\h{F}_r=F_r-\frac{1}{2}\{\theta\wedge\theta_r+\theta_s\wedge\theta_t\}$ implies that
\begin{equation}
\h{F}_2+i\h{F}_3=F_2+iF_3-\frac{1}{2}(\theta+i\theta_1)\wedge(\theta_2+i\theta_3).
\end{equation}
It shows that the function $f_i$ satisfies the condition in (\ref{pott}). 
\eproof

Next, we investigate the geometry of the Lee field with respect to the geometry of
the HKT-metric ${\hat g}$ and its associated HKT-connection $D$. The following
result can be found in \cite{PPS1}.
\begin{proposition}\label{symmetry}
 Let $V$ be the vector field dual to the parallel Lee-form
with respect to the locally conformally hyperk\"ahler metric $g$, then the algebra
$\{V\}\oplus\{I_1V, I_2V, I_3V\}$ is isomorphic to $\lie{u}(1)\oplus\lie{su}(2)$. 
Moreover,
\begin{equation}
{\mathcal L}_VI_r=0, \quad {\mathcal L}_Vg=0,\quad {\mathcal L}_{I_rV}g=0, \quad {\mathcal L}_{I_rV}I_s=\epsilon^{rst}I_t
.
\end{equation}
\end{proposition}

To understand the relation between HKT-geometry and the Lee field $V$, 
we need to describe the behavior of the Lee field with respect to the forms $\theta$ and $\theta_r$.

\begin{lemma} Let $V$ be the Lee field, $\theta_r=I_r\theta$ for $1\leq r\leq 3$. Then
\begin{equation}
{\mathcal L}_V\theta=0, \quad {\mathcal L}_V\theta_r=0,
\quad
{\mathcal L}_{I_rV}\theta=0, \quad {\mathcal L}_{I_rV}\theta_s=\epsilon^{rst}\theta_t.
\end{equation}
\end{lemma}
\bproof
The Lee form $\theta$ is invariant along its dual vector field because it is parallel with respect to 
the Levi-Civita connection of the locally conformally hyperk\"ahler metric $g$. The forms $\theta_r$ are invariant with respect to the
Lee field because the Lee form is invariant and the Lee field is hypercomplex.

Next, for any vector field $Y$,
\begin{eqnarray*}
({\mathcal L}_{I_rV}\theta)Y &=& I_rV(\theta(V))-\theta({\mathcal L}_{I_rV}Y)=I_rVg(V, Y)-g(V, [I_rV, Y])\\
&=& g(\nabla_{I_rV}V, Y)+g(V, \nabla_{I_rV}Y)-g(V, [I_rV, Y])=g(V, \nabla_{I_rV}Y-[I_rV, Y])\\
&=&g(V, \nabla_Y(I_rV))=Yg(V, I_rV)-g(\nabla_YV, I_rV)=0.
\end{eqnarray*}
It follows that ${\mathcal L}_{I_rV}\theta=0$. This equality is combined with ${\mathcal L}_{I_rV}I_s=\epsilon^{rst}I_t$
to yield the last one in this lemma.
\eproof

Due to Lemma \ref{parallelism}, we learn the following.
\begin{theorem}\label{D-parallel} The potential 1-form for the HKT-metric ${\hat g}$ is
parallel. 
\end{theorem}
\bproof
The tensor $\theta^2+\theta_1^2+\theta_2^2+\theta_3^2$
is invariant with respect to the given vector fields due to the last lemma. As
${\mathcal L}_Vg=0$ and ${\mathcal L}_{I_rV}g=0$,
the vector fields $V, I_1V, I_2V, I_3V$ are Killing vector fields of the HKT-metric
$\hat g$. By Lemma \ref{parallelism}, the vector field $V$ is parallel with
respect to the HKT-connection $D$. Since $D$ is a Riemannian connection, 
the dual 1-form $\htheta$ is parallel.
\eproof

\subsection{Additional examples of HKT-spaces with parallel potential 1-form}\label{product}
Once we construct HKT-spaces with $D(2,1;-1)$-symmetry, we 
 can
generate new examples
through direct 
 products.
Indeed let \( (M_1,g_1,I_r^{(1)}) \), \( (M_2,g_2,I_r^{(2)}) \) be two
locally conformally hyperk\"ahler manifolds with parallel Lee forms.  Then
\( \hat g_i \) are HKT metrics with special homotheties \( V_i \), \( i=1,2
\).  On \( M=M_1\times M_2 \) consider the product metric 
\begin{equation}
  \hat g = \frac12(\pi_1^*\hat g_1 + \pi_2^*\hat g_2)
\end{equation}
and complex structures 
$  I_r = (I_r^{(1)},I_r^{(2)})$.
This geometry on \( M \) is HKT, since 
\[  F_r = \frac12(\pi_1^*F^{(1)}_r + \pi_2^*F^{(2)}_r)
\]
and 
$  c= -d_rF_r = -I_rdF_r = \frac12(\pi_1^*c_1+\pi_2^*c_2)$
is independent of~\( r = 1,2,3 \).
Let 
\begin{equation}
  V = (V_1, V_2),\qquad \htheta = \frac12(\pi_1^*\htheta^{(1)}+\pi_2^*\htheta^{(2)}).
\end{equation}
Then \( V \) generates a $D(2,1;-1)$-symmetry, since this is true
of \( V_1 \) and \( V_2 \).  Moreover, $\htheta$ is a potential 1-form.
 Note that the normalization of $\hat g$ has been chosen to fit with
 conventions of the next section.

\section{Relating Torsion 3-Forms and Potential 1-Forms}

The past section demonstrates that locally conformally hyperk\"ahler
manifolds with parallel Lee form generate HKT-spaces with
$D(2,1;-1)$-symmetry.  In this section, we demonstrate that the latter type
of geometry is more general than the former. This is achieved through an
analysis of 
 the
torsion 3-form.

Consider now an HKT structure obtained from a locally conformally
hyperk\"ahler metric with parallel Lee form.  The torsion three-form is
given by the following lemma.

\begin{lemma}\label{cubic form}
  The torsion three-form is determined by \( \hat \theta \) as
  \begin{equation}
    \label{eq:c-hat}
    c = -(\hat\theta_1\wedge \hat F_1 + \hat\theta_2\wedge \hat F_2 +
    \hat\theta_3\wedge \hat F_3 -
    2\hat\theta_1\wedge\hat\theta_2\wedge\hat\theta_3)
  \end{equation}
\end{lemma}
\bproof 
To calculate the torsion 3-form when the HKT-structure is generated by a locally conformally
hyperk\"ahler structure,
we recall \( \hat\theta=\frac12\theta \). Next, we write
equation~\eqref{eq:hat-F} as
\begin{equation}
  F_1 = \hat F_1 +
  2(\hat\theta\wedge\hat\theta_1+\hat\theta_2\wedge\hat\theta_3). 
\end{equation}
Then from equation~\eqref{eq:c-theta-F}, we have
\begin{eqnarray*}
  c
  &=&-I_1d{\hat F}_1
   = -\frac12(\theta_1\wedge F_1 + \theta_2\wedge F_2 + \theta_3\wedge F_3
  - 2 \theta_1 \wedge \theta_2
  \wedge\theta_3 )
  \\
  &=&-(\hat\theta_1\wedge \hat F_1 + \hat\theta_2\wedge \hat F_2 +
    \hat\theta_3\wedge \hat F_3 -
    2\hat\theta_1\wedge\hat\theta_2\wedge\hat\theta_3),
  \end{eqnarray*}
  as claimed.
Thus the torsion is an inhomogeneous cubic function of the one-form \(
\htheta \).  
\eproof

The torsion
three-form \( c \) determines a torsion one-form \( \tau \) by 
\begin{equation}
  \tau(X) = \frac12 \sum_{i=1}^{4m} c(I_r X,e_i,I_r e_i).
\end{equation}
The HKT condition ensures that \( \tau \) is independent of the choice of
\( I_r \), \( r=1,2,3 \) \cite{I}.  Under the current constraints, 
\begin{equation}
\tau(X)=(2m-1+\norm\htheta^2)\htheta(X) \label{eq:htheta-tau}.
\end{equation}
Thus \( \htheta = \lambda\tau \), where \( \lambda \) is the unique real
(and positive) solution to the cubic equation
\begin{equation}
  \lambda(2m-1+\lambda^2) = 1.
\end{equation}
On an arbitrary HKT manifold, whose torsion one-form is non-zero, one may
always find a one-form \( \htheta \) satisfying~\eqref{eq:htheta-tau}.  By
rescaling \( \hat g \) by a homothety, we may ensure that \(
\norm{\htheta}^2 = 1/2 \) at some base point.  With these conventions we
call \( \htheta \) a \emph{normalized torsion one-form} of \( M \).  We
say that an HKT manifold \( M \) is of \emph{cubic type} if its torsion
three-form~\( c \) is related to the normalized torsion one-form \( \htheta
\) by equation~\eqref{eq:c-hat}.

Let \( V \) be the vector field dual to \( \htheta \) via \( \hat g \) in
this normalization.  Then \( \htheta = \hat g(V,\cdot) \) and
\begin{equation}
  \label{eq:length}
  \hat g(V,V)=\frac12,\quad\mbox{or equivalently},\quad \htheta(V)=\frac12.
\end{equation}

\begin{theorem}\label{inverse}
  Suppose \( (M,\hat g,I_1,I_2,I_3) \) is an HKT manifold with a normalized
  torsion one-form \( \htheta \).  If the torsion~\( c \) is given by
  \begin{equation}
    c = -\{\hat\theta_1\wedge \hat F_1 + \hat\theta_2\wedge \hat F_2 +
    \hat\theta_3\wedge \hat F_3 -
    2\hat\theta_1\wedge\hat\theta_2\wedge\hat\theta_3\}
  \end{equation}
  and the dual vector field of the torsion 1-form generates a
 $D(2,1;-1)$-symmetry,  then
  \begin{equation}\label{sheared}
    g = \hat g + 2\{\htheta\otimes \htheta + \htheta_1\otimes
    \htheta_1 + \htheta_2\otimes \htheta_2 + \htheta_3\otimes
    \htheta_3 \}
  \end{equation}
  is locally conformally hyperk\"ahler with parallel Lee form.
\end{theorem}
\bproof We first compute the derivatives of \( \htheta \) and \( \htheta_r
\). Let $V$ be the dual vector field of the 1-form $\htheta$. By definition of symmetry
and Lemma \ref{parallelism}, \( V \) is parallel. By Lemma~\ref{diff}, we have
\begin{equation}
  d\htheta(X,Y)=c(V,X,Y),\qquad\htheta_1 (X,Y) = c(I_1V,X,Y).
\end{equation}
The form of \( c \) now gives
\begin{eqnarray*}
    d\htheta_1
  &=& -\left(\frac12\hat F_1 - \htheta_1\wedge F_1(I_1V,\cdot)
  -\htheta_2\wedge F_2(I_1V,\cdot) - 
  \htheta_3\wedge F_3(I_1V,\cdot) -\htheta_2\wedge\htheta_3\right)\\
  &=& -\frac12\hat F_1 + \htheta\wedge\htheta_1 -
  \htheta_2\wedge\htheta_3 
  = -\frac12 F_1 + \frac14\{\theta\wedge\theta_1 +
  \theta_2\wedge\theta_3\} - \frac14\theta\wedge\theta_1 +
  \frac14\theta_2\wedge\theta_3\\
  &=& -\frac12 F_1 + \frac12\theta\wedge\theta_1,
\end{eqnarray*}
where \( \theta=2\htheta \) and \( F_1= g(I_1\cdot,\cdot) \) is given
by~\eqref{eq:hat-F}.  Thus \( F_1 = \theta\wedge\theta_1-d\theta_1 \) and
this has derivative
\begin{equation}
  dF_1 = d(\theta\wedge\theta_1-d\theta_1) = -\theta\wedge d\theta_1 =
  \theta\wedge F_1.
\end{equation}
As similar equations hold for \( F_2 \) and \( F_3 \), we conclude that \(
g \) is locally conformally hyperk\"ahler.  The Lee form is a constant
multiple of \( \theta \), which is closed and hence parallel.
\eproof

The condition on the structure of the torsion three-form 
 is 
rather strong.
However, this is a necessary condition. The example in Section \ref{product} 
demonstrates that the existence of $D(2,1;-1)$-symmetry itself does not 
necessarily come from a locally conformally hyperk\"ahler manifold. This is 
consistent with the fact that in general the product of locally conformally K\"ahler 
manifolds is not necessarily locally conformally K\"ahler. 
In fact,
 the torsion of the example given in Section \ref{product}
is 
not of cubic type.
If we consider 
 the case where each factor is locally conformally hyperk\"ahler, put
$  g = \hat g + 2\{\htheta\otimes \htheta + \htheta_1\otimes
  \htheta_1 + \htheta_2\otimes \htheta_2 + \htheta_3\otimes
  \htheta_3 \}
$
and \( \theta=2\htheta \), the K\"ahler form 
 $ F_1$ is equal to 
\[
\frac12\left(\pi_1^*F_1^{(1)} + \pi_2^*F_2^{(2)} + 
  \pi_1^*\theta^{(1)}\wedge
  \pi_2^*\theta_1^{(2)}+\pi_1^*\theta_2^{(1)}\wedge\pi_2^*\theta_3^{(2)}
  + \pi_2^*\theta^{(2)}\wedge
  \pi_1^*\theta_1^{(1)}+\pi_2^*\theta_2^{(2)}\wedge\pi_1^*\theta_3^{(1)}
  \right),
\]
so 
\begin{eqnarray*}
 2 dF_1
  &=& \pi_1^*\theta^{(1)}\wedge F_1^{(1)} + \pi_2^*\theta^{(2)}\wedge
  F_1^{(2)}\\
  &&\qquad
  -\pi_1^*\theta^{(1)}\wedge
  \pi_2^*(\theta^{(2)}\wedge\theta_1^{(2)}-F_1^{(2)})
  +\pi_1^*(\theta^{(1)}\wedge\theta_2^{(1)}-F_2^{(1)})
  \wedge\pi_2^*\theta_3^{(2)}\\  
  &&\qquad
  -\pi_1^*\theta_2^{(1)} \wedge
  \pi_2^*(\theta^{(2)}\wedge\theta_3^{(2)}-F_3^{(2)})
  -\pi_2^*\theta^{(2)}\wedge
  \pi_1^*(\theta^{(1)}\wedge\theta_1^{(1)}-F_1^{(1)})\\
  &&\qquad
  +\pi_2^*(\theta^{(2)}\wedge\theta_2^{(2)}-F_2^{(2)})
  \wedge\pi_1^*\theta_3^{(1)}
  -\pi_2^*\theta_2^{(2)} \wedge
  \pi_1^*(\theta^{(1)}\wedge\theta_3^{(1)}-F_3^{(1)})
  \\
  &\ne& 2\theta\wedge F_1,
\end{eqnarray*}
since the expression contains non-zero terms involving for example \(
\pi_1^*F_2^{(1)} \)
 and terms such as \( \theta^{(1)}\wedge F_1^{(2)} \)
occur with the wrong coefficients
.
Thus \( g \) is \emph{not} locally conformally hyperk\"ahler.

\bigbreak\noindent{\bf Remark:} 
  There is an alternative way to see when an HKT-space with
  $D(2,1;-1)$-symmetry will generate a locally conformally hyperk\"ahler
  metric using the transformation of the last theorem. Suppose that the
  dual vector field of a closed 1-form $\htheta$ is a $D(2,1,;-1)$-symmetry
  on an HKT-space. Now we do not assume that the torsion of the HKT-space
  is of cubic type.  Define $\theta=2{\hat\theta}$. By Proposition
  \ref{converse}, $-\theta$ is a potential 1-form for the HKT metric $\hat
  g$.  Again, consider the Riemannian metric (\ref{sheared}).  Due to the
  choice of $V$, $\theta$ is the dual of the vector field $V$ with respect
  to the metric $g$.  Define
  $g_0=\theta\otimes\theta+\theta_1\otimes\theta_1+\theta_2\otimes\theta_2
  +\theta_3\otimes\theta_3.$ Then for any vector fields $X$ and $Y$, when
  $rst$ is a cyclic permutation of $123$,
\[
g_0(I_rX, Y)=(\theta\wedge\theta_r+\theta_s\wedge\theta_t)(X,Y).
\]
Therefore, $ F_r={\hat
F}_r+\frac12(\theta\wedge\theta_r+\theta_s\wedge\theta_t) ={\hat
F}_r+2(\htheta\wedge\htheta_r+\htheta_s\wedge\htheta_t)$. Since $-\theta$
is a potential 1-form,
\begin{equation}
{\hat F}_r =-\frac12 (d\theta_r+d_s\theta_t)=-\frac12(d\theta_r-I_sd\theta_r)
     =-\frac12(d\theta_r-I_td\theta_r).
\end{equation}
It follows that
\begin{eqnarray*}
F_r&=&-\frac12(d\theta_r-I_sd\theta_r)
 +\frac12 (\theta\wedge\theta_r+\theta_s\wedge\theta_t)
= -\frac12\{(d\theta_r-\theta\wedge\theta_r)
      -I_s(d\theta_r-\theta\wedge\theta_r)\}\\
\mbox{and } &=& -\frac12\{(d\theta_r-\theta\wedge\theta_r)
      -I_t(d\theta_r-\theta\wedge\theta_r)\}.
\end{eqnarray*}
Therefore, $F_r=-(d\theta_r-\theta\wedge\theta_r)$ if and only if for
$s\neq r$, $
I_s(d\theta_r-\theta\wedge\theta_r)=-(d\theta_r-\theta\wedge\theta_r).  $
On the other hand, we check that $ I_a(d\theta_a-\theta\wedge\theta_a) =
d\theta_a-\theta\wedge\theta_a.  $ The conclusion is the following
observation.

\begin{proposition} The metric
$g$ is a locally conformal hyperk\"ahler
 metric with parallel Lee form $\theta$ if and only if for all
$s\neq r$,
$I_s(d\theta_r-\theta\wedge\theta_r)=-(d\theta_r-\theta\wedge\theta_r).$
\end{proposition}

\noindent{\bf Remark:} 
An HKT-structure is said to be strong if the torsion 3-form $c$ is closed
\cite{GPS} \cite{MS}.  We calculate exterior differential of the torsion
3-form when the HKT-structure is generated by a locally conformally
hyperk\"ahler structure. We continue to use the notation in Lemma
\ref{cubic form}. With the aid of (\ref{dFr}) and (\ref{trei}),
\begin{eqnarray*}
dc &=& -\frac12 
 (d\theta_1\wedge F_1+d\theta_2\wedge F_2+d\theta_3\wedge F_3
    -\theta_1\wedge dF_1-\theta_2\wedge dF_2-\theta_3\wedge dF_3
  \nonumber\\
 & & -2d\theta_1\wedge \theta_2\wedge\theta_3
    +2\theta_1\wedge d\theta_2\wedge\theta_3
     -2\theta_1\wedge \theta_2\wedge d\theta_3)
\nonumber\\
     & =& \frac12 \left(
     \left(F_1-\theta\wedge\theta_1-\theta_2\wedge\theta_3\right)^2
     +\left(F_2-\theta\wedge\theta_2-\theta_3\wedge\theta_1\right)^2
     +\left(F_3- \theta\wedge\theta_3-\theta_1\wedge\theta_2\right)^2
   \right).
  \end{eqnarray*}
  
  This formula demonstrates that the restriction of $dc$ on the
  quaternionic span of $V$ is equal to zero.  On the quaternionic
  complement it is equal to
\begin{equation}
\frac12 (F_1\wedge F_1+F_2\wedge F_2+F_3\wedge F_3).
\end{equation}
In particular, it shows the following observation.

\begin{proposition} If $M$ is a locally conformally hyperk\"ahler space
  with real dimensional at least 8, then the associated HKT-structure $\hat
  g$ is never strong.
\end{proposition}

\end{document}